\documentclass[11pt, reqno, psamsfonts]{amsart}
\pdfoutput=1

\usepackage{amssymb}
\usepackage{amsthm}
\usepackage{amsmath}
\usepackage{latexsym}
\usepackage[T1]{fontenc}
\usepackage[utf8]{inputenc}
\usepackage[russian, french, english]{babel}

\usepackage{graphicx}
\usepackage{wrapfig}
\usepackage[justification=centering, labelfont=bf]{caption}
\usepackage{mathtools}
\usepackage[hidelinks]{hyperref}
\usepackage{tikz}
\usepackage{amsbsy}
\usepackage[inline]{enumitem}
\usepackage{mathrsfs}
\usetikzlibrary{shapes,snakes}
\usetikzlibrary{arrows.meta}
\usetikzlibrary{decorations.pathreplacing}
\usepackage{float}
\usepackage{multicol}
\usepackage{bbm}
\usepackage{bm}
\usepackage{stmaryrd}
\usepackage{mathabx}\usepackage{lmodern}
\usepackage{titlesec}

\usepackage[spacing=true,kerning=true,babel=true,tracking=true]{microtype}

\usepackage[shortcuts]{extdash}

\usepackage[foot]{amsaddr} 

\usepackage[left=1in,right=1in,top=1in,bottom=1in,bindingoffset=0cm]{geometry}

\usepackage[backend=biber, style=alphabetic, sorting=nyt, maxnames=100]{biblatex}
\title{A Short Proof of Bernoulli Disjointness via the\\ Local Lemma}
\date{}
\author{Anton~Bernshteyn}
\address{\normalfont{}\textls{Department of Mathematical Sciences, Carnegie Mellon University, Pittsburgh, PA, 15213, USA}}
\email{abernsht@math.cmu.edu}

\newtheoremstyle{bfnote}%
{}{}%
{\slshape}{}%
{\bfseries}{\bfseries.}%
{ }%
{\thmname{#1}\thmnumber{ #2}\thmnote{ \ep{\normalfont{}#3}}}

\newtheoremstyle{defbfnote}%
{}{}%
{}{}%
{\bfseries}{.}%
{ }%
{\thmname{#1}\thmnumber{ #2}\thmnote{ (#3)}}

\newtheoremstyle{claim}%
{}{}%
{\slshape}{}%
{\itshape}{.}%
{ }%
{\thmname{#1}\thmnumber{ #2}\thmnote{ \ep{\normalfont{}#3}}}

\theoremstyle{bfnote}

\newtheorem{theo}[equation]{Theorem}

\newtheorem{lemma}[equation]{Lemma}

\newtheorem*{claim*}{Claim}
\newtheorem*{corl*}{Corollary}

\theoremstyle{claim}

\newcommand*{\myproofname}{Proof}

\theoremstyle{definition}

\newtheorem*{defn*}{Definition}

\newtheorem*{exmp*}{Example}

\theoremstyle{remark}
\newtheorem*{ques*}{Question}
\newtheorem*{remk*}{Remark}

\newcommand{\0}{\varnothing}
\newcommand{\set}[1]{\{#1\}}

\newcommand{\dom}{\mathrm{dom}}

\newcommand{\acts}{\curvearrowright}
\newcommand{\N}{{\mathbb{N}}}
\newcommand{\Z}{\mathbb{Z}}
\newcommand{\F}{\mathbb{F}}
\newcommand{\R}{\mathbb{R}}
\renewcommand{\C}{\mathbb{C}}
\newcommand{\Q}{\mathbb{Q}}

\renewcommand{\epsilon}{\varepsilon}
\renewcommand{\phi}{\varphi}
\renewcommand{\theta}{\vartheta}
\newcommand{\E}{\mathbb{E}}
\renewcommand{\leq}{\leqslant}
\renewcommand{\geq}{\geqslant}

\renewcommand{\G}{\Gamma}

\newcommand{\defeq}{\coloneqq}

\newcommand{\emphd}[1]{{\fontseries{b}\selectfont\textsf{#1}}}
\newcommand{\B}{{\mathscr{B}}}
\renewcommand{\P}{\mathbb{P}}

\newcommand{\Nbhd}{\mathrm{N}}

\newcommand{\bemph}[1]{{\normalfont#1}} 
\newcommand{\ep}[1]{\bemph{(}#1\bemph{)}} 

\newenvironment{scproof}[1][Proof]{\begin{proof}[\textsc{#1}]}{\end{proof}}

\usepackage{xspace}

\numberwithin{equation}{section}

\renewcommand{\thesubsection}{\arabic{section}.\Alph{subsection}}

\usepackage{etoolbox}

\titleformat{\section}[block]{\scshape\filcenter}{\thesection.}{1ex}{}
\titleformat{\subsection}[block]{\bfseries\filcenter}{\thesubsection.}{1ex}{}
\titleformat{\subsubsection}[runin]{\bfseries}{\bfseries\upshape\thesubsubsection.}{1ex}{}[---]

\titlespacing*{\section}{0pt}{*3}{*1}
\titlespacing*{\subsection}{0pt}{*3}{*1}
\titlespacing*{\subsubsection}{0pt}{*1.5}{*0}

\makeatletter
\newcommand{\neutralize}[1]{\expandafter\let\csname c@#1\endcsname\count@}
\makeatother

\renewcommand{\mkbibnamefirst}[1]{\textsc{#1}}
\renewcommand{\mkbibnamelast}[1]{\textsc{#1}}
\renewcommand{\mkbibnameprefix}[1]{\textsc{#1}}
\renewcommand{\mkbibnameaffix}[1]{\textsc{#1}}
\renewcommand{\finentrypunct}{}

\renewbibmacro{in:}{}

\renewbibmacro*{volume+number+eid}{%
	\printfield{volume}%
	\setunit*{\addnbspace}
	\printfield{number}%
	\setunit{\addcomma\space}%
	\printfield{eid}}

\DeclareFieldFormat[article]{volume}{\textbf{#1}\space}
\DeclareFieldFormat[article]{number}{\mkbibparens{#1}}

\DeclareFieldFormat{journaltitle}{#1,}
\DeclareFieldFormat[thesis]{title}{\mkbibemph{#1}\addperiod}
\DeclareFieldFormat[article, unpublished, thesis]{title}{\mkbibemph{#1},}
\DeclareFieldFormat[book]{title}{\mkbibemph{#1}\addperiod}
\DeclareFieldFormat[unpublished]{howpublished}{#1, }

\DeclareFieldFormat{pages}{#1}

\DeclareFieldFormat[article]{series}{Ser.~#1\addcomma}

\renewcommand*{\bibfont}{\small}

\newcommand\blfootnote[1]{%
	\begingroup
	\renewcommand\thefootnote{}\footnote{#1}%
	\addtocounter{footnote}{-1}%
	\endgroup
}

\begin{document}
	
	\pagestyle{plain}
	
	\maketitle
	
	\blfootnote{2010 Mathematics Subject Classification: Primary 37B05, 37B10; Secondary 05D40.}%
	
	\blfootnote{Key words and phrases: disjointness, minimal flows, Bernoulli flow, Lov\'asz Local Lemma.}%
	
	\begin{abstract}
		Recently, Glasner, Tsankov, Weiss, and Zucker showed that if $\G$ is an infinite discrete group, then every minimal $\G$-flow is disjoint from the Bernoulli shift $2^\G$. Their proof is somewhat involved; in particular, it invokes separate arguments for different classes of groups. In this note, we give a short and self-contained proof of their result using purely combinatorial methods  applicable to all groups at once. Our proof relies on the Lov\'asz Local Lemma, an important tool in probabilistic combinatorics that has recently found several applications in the study of dynamical systems.
	\end{abstract}

	\section{Introduction}
	
	Throughout, $\G$ denotes an infinite discrete \ep{but not necessarily countable} group. All topological spaces in this paper are assumed to be Hausdorff. A \emphd{$\G$-flow} is a nonempty compact space $X$ equipped with a continuous action $\G \acts X$. A \emphd{subflow} of a $\G$-flow $X$ is a nonempty closed $\G$-invariant subset $Y \subseteq X$. A $\G$-flow $X$ is \emphd{minimal} if there is no subflow $Y \subseteq X$ with $Y \neq X$; equivalently, $X$ is minimal if the orbit of every point $x \in X$ is dense in $X$. A straightforward application of Zorn's lemma shows that every $\G$-flow has a minimal subflow. This demonstrates that minimal flows \emph{exist}, but in general it is very difficult to say much about their structure and behavior.
	
	Let $X$ be a nonempty compact space. The space $X^\G$ of all functions $\G\to X$, equipped with the product topology, is also compact, and it becomes a $\G$-flow under the action $\G \acts X^\G$ given by
	\[
		(\gamma \cdot x)(\delta) \,\defeq\, x(\delta\gamma) \qquad \text{for all $\gamma$, $\delta \in \G$ and $x \in X^\G$}.
	\]
	The $\G$-flows of the form $X^\G$ are called \emphd{Bernoulli shifts}, or simply \emphd{shifts}. A particularly important case is when $X$ is a finite set with the discrete topology; for concreteness, we may then assume that $X = k$ for some $k \in \N$. \ep{Here we identify each $k \in \N$ with the $k$-element set $\set{i \in \N \,:\, i < k}$.}
	
	The following notions were introduced by Furstenberg in \cite{Furstenberg}. Let $X$ and $Y$ be $\G$-flows. We view the product space $X \times Y$ as a $\G$-flow equipped with the diagonal action of $\G$. A \emphd{joining} of $X$ and $Y$ is a subflow $Z \subseteq X \times Y$ that projects onto $X$ and $Y$. The $\G$-flows $X$ and $Y$ are \emphd{disjoint}, in symbols $X \perp Y$, if they have only one joining, namely $X \times Y$ itself. It is not hard to see that if $X \perp Y$, then at least one of $X$ and $Y$ is minimal.
	
	Recently, Glasner, Tsankov, Weiss, and Zucker obtained the following result:

	\begin{theo}[{Glasner--Tsankov--Weiss--Zucker \cite{GTWZ}}]\label{theo:disj}
		If $X$ is a minimal $\G$-flow, then $X \perp 2^\G$.
	\end{theo}

	For $\G = \Z$, Theorem~\ref{theo:disj} was proven earlier by Furstenberg \cite{Furstenberg}. 
	
	The first step in the proof of Theorem~\ref{theo:disj} in \cite{GTWZ} is a reduction to a combinatorial question concerning the so-called \emph{separated covering property} of minimal $\G$-flows. Somewhat surprisingly, to answer this question, \cite{GTWZ} makes heavy use of group theory; in particular, two classes of groups \ep{namely \emph{maximally almost periodic} and \emph{ICC} groups} are treated separately and with completely different arguments. The proof in the ICC case relies on the recent breakthrough results of Frisch, Tamuz, and Vahidi Ferdowsi \cite{FTF}.
	
	In this note, we give a short and self-contained proof of Theorem~\ref{theo:disj} that is purely combinatorial and treats all groups simultaneously. We deduce Theorem~\ref{theo:disj} from the following fact, which is interesting in its own right. Let $X$ be a $\G$-flow. We say that a set $U \subseteq X$ \emphd{traps} a point $x \in X$ \ep{or that $x$ is \emphd{trapped} in $U$} if the orbit of $x$ is contained in $U$.

	\begin{theo}\label{theo:synd}
		Let $k \in \N^+$ and let $V \subseteq k^\G$ be a nonempty open set. Then there exists $n \in \N$ such that for all finite $F \subset \G$ of size at least $n$, the set $F \cdot V$ traps a point.
	\end{theo}

	Theorem~\ref{theo:disj} is derived from Theorem~\ref{theo:synd} in \S\ref{sec:reduction}. To establish Theorem~\ref{theo:synd}, we rely on the so-called \emph{Lov\'asz Local Lemma} \ep{the \emph{LLL} for short}, an important tool in probabilistic combinatorics introduced by Erd\H{o}s and Lov\'asz in \cite{EL}.
	We state the LLL, in the form we will need, in \S\ref{sec:LLL} and use it to prove Theorem~\ref{theo:synd} in \S\ref{sec:synd}. While the LLL has been widely used in combinatorics and graph theory for over forty years now, it has only recently become apparent that the LLL can be applied to the study of dynamical systems as well. Several applications of the LLL in ergodic theory and topological dynamics can be found in \cite{ABT, LSS, MLLL, AW}; this paper contributes yet another one.

	\section{Derivation of Theorem~\ref{theo:disj} from Theorem~\ref{theo:synd}}\label{sec:reduction}
	
	Let $X$ be a minimal $\G$-flow and let $Z \subseteq X \times 2^\G$ be a joining. Our goal is to show that $Z = X \times 2^\G$. To this end, let $U \subseteq X$ and $V \subseteq 2^\G$ be nonempty open sets. We have to argue that $Z \cap (U \times V) \neq \0$.
	
	\begin{lemma}\label{lemma:rec}
		If $X$ is a minimal $\G$-flow and $U \subseteq X$ is a nonempty open set, then there exists an infinite subset $S \subseteq \G$ with $\bigcap \set{\sigma \cdot U \,:\, \sigma \in S} \neq \0$.
	\end{lemma}
	\begin{scproof}
		Since $X$ is minimal, the orbit of each point in $X$ intersects $U$, and thus $X = \G \cdot U$. The compactness of $X$ implies that there is a finite subset $T \subset \G$ such that $X = T \cdot U$. Say that a set $L \subseteq \G$ is \emphd{right $T$-separated} if $\lambda_1T \cap \lambda_2T = \0$ for every pair of distinct $\lambda_1$, $\lambda_2 \in L$. Since $T$ is finite while $\G$ is infinite, there is an infinite right $T$-separated subset $L \subseteq \G$. Consider any $x \in U$. For each $\lambda \in L$, we have $\lambda^{-1} \cdot x \in T \cdot U$, i.e., $x \in \lambda T \cdot U$. In other words, for each $\lambda \in L$, there is $\tau_\lambda \in T$ with $x \in \lambda \tau_\lambda \cdot U$. Since $L$ is right $T$-separated, the set $S \defeq \set{\lambda \tau_\lambda \,:\, \lambda \in L}$ is infinite, and we are done. 
	\end{scproof}
	
	Let $S \subseteq \G$ be an infinite set as in Lemma~\ref{lemma:rec} applied to $U$, and let $n \in \N$ be the quantity given by Theorem~\ref{theo:synd} applied to $V$. Pick an arbitrary finite subset $F \subset S$ of size at least $n$ and let
	\[
	U^\ast \,\defeq\, \bigcap \set{\sigma \cdot U \,:\, \sigma \in F}.
	\]
	By the choice of $n$, the set $F \cdot V$ traps a point, say $y$. Since $Z$ projects onto $2^\G$, there is some $x \in X$ with $(x,y) \in Z$. As $X$ is minimal and $U^\ast$ is nonempty open, there is $\gamma \in \G$ with $\gamma \cdot x \in U^\ast$. Since $y$ is trapped in $F \cdot V$, we have $\gamma \cdot y \in F \cdot V$, which means that $\sigma^{-1}\gamma \cdot y \in V$ for some $\sigma \in F$, and since $\gamma \cdot x \in U^\ast$, we also have $\sigma^{-1}\gamma \cdot x \in U$. Therefore,
	\[
	(\sigma^{-1}\gamma \cdot x, \, \sigma^{-1}\gamma \cdot y) \,=\, \sigma^{-1}\gamma \cdot (x,y) \,\in\, Z \cap (U \times V),
	\]
	and the proof of Theorem~\ref{theo:disj} is complete.
	
	\section{The Lov\'asz Local Lemma}\label{sec:LLL}
	
	We shall only require a somewhat specialized but simplified version of the Lov\'asz Local Lemma; for a more general discussion of the LLL, the reader is referred to~\cite{AS} and \cite{MR}. The presentation below follows, with slight modifications,~\cite[\S{}1.2]{MLLL}.
	
	Let $X$ be an arbitrary set and let $k \in \N^+$. A \emphd{bad \ep{$k$-}event} over $X$ is a nonempty set $B$ of partial functions $\phi \colon X \rightharpoonup k$ with finite domains such that for all $\phi$, $\phi' \in B$, $\dom(\phi) = \dom(\phi')$. If a bad event $B$ is nonempty, then its \emphd{domain} is the set $\dom(B) \coloneqq \dom(\phi)$ for any \ep{hence all} $\phi \in B$; the domain of the empty bad event is, by definition, the empty set. The \emphd{probability} of a bad $k$-event $B$ with domain $D$ is defined to be $\P[B] \coloneqq |B|/k^{|D|}$. A map $f \colon X \to k$ \emphd{avoids} a bad event $B$ if there is no $\phi \in B$ such that $\phi \subseteq f$. Notice that if $X$ is finite and $f \colon X \to k$ is drawn uniformly at random from $k^X$, then $\P[B]$ is precisely the probability that $f$ \emph{does not} avoid $B$.
	
	A~\emphd{\ep{$k$\=/}instance \ep{of the LLL}} over a set $X$ is an arbitrary set $\B$ of bad $k$-events over $X$. A~\emphd{solution} to a $k$-instance~$\B$ is a function $f \colon X \to k$ that avoids every $B \in \B$. For an instance $\B$ and $B \in \B$, the \emphd{neighborhood} of $B$ in $\B$ is the set
	\[
	\Nbhd_\B(B) \defeq \set{B' \in \B \setminus \set{B} \,:\, \dom(B') \cap \dom(B) \neq \0}.
	\]
	The \emphd{degree} of $B$ in $\B$ is defined to be $\deg_\B(B) \defeq |\Nbhd_\B(B)|$. Let
	\[
	p(\B) \defeq \sup_{B \in \B} \P[B] \qquad \text{and} \qquad d(\B) \defeq \sup_{B \in \B} \deg_\B(B).
	\]
	An instance $\B$ is \emphd{correct for the LLL}, or simply \emphd{correct}, if
	\[
	e\cdot p(\B) \cdot (d(\B) + 1) < 1,
	\]
	where $e = 2.71\ldots${} denotes the base of the natural logarithm.
	
	\begin{theo}[{Erd\H os--Lov\'asz~\cite{EL}; \textls{Lov\'asz Local Lemma}}]\label{theo:LLL}
		Let $k \in \N^+$ and let $\B$ be a $k$-instance of the LLL over a set $X$. If $\B$ is correct for the LLL, then $\B$ has a solution.
	\end{theo}
	
	The LLL was introduced by Erd\H os and Lov\'asz \ep{with $4$ in place of $e$} in their seminal paper~\cite{EL}; the constant was subsequently improved by Lov\'asz (the sharpened version first appeared in~\cite{Spencer}).
	
	We should mention that there are two respects in which Theorem~\ref{theo:LLL} in the above form is less general then the ``full'' LLL. First, Theorem~\ref{theo:LLL} only works with \emph{product} probability spaces such as $k^X$; this is a special case of the~LLL in the so-called \emph{variable framework} (the~name is due to Kolipaka and Szegedy~\cite{KolipakaSzegedy}).
	However, although this case is special, it does encompass most typical applications. For the statement of the~LLL for general probability spaces, see~\cite[Corollary~5.1.2]{AS}. Deducing Theorem~\ref{theo:LLL} from \cite[Corollary~5.1.2]{AS} is routine when $X$ is finite (see, e.g., \cite[41]{MR}); the case of infinite~$X$ then follows by compactness.
	
	Second, there is a more general form of the LLL \ep{often referred to as the \emph{General Lov\'asz Local Lemma}}, that applies to instances without a uniform upper bound on $\deg_\B(B)$; see \cite[Theorem 5.1.1]{AS}. However, this more general statement is somewhat technical and we will not need it here.
	
	\section{Proof of Theorem~\ref{theo:synd}}\label{sec:synd}
	
	Let $D \subset \G$ be a finite set. We say that a set $L \subseteq \G$ is \emphd{left $D$-separated} if $D\lambda_1 \cap D\lambda_2 = \0$ for every pair of distinct $\lambda_1$, $\lambda_2 \in L$.
	
	\begin{lemma}\label{lemma:indep}
		If $D$, $F \subset \G$ are finite sets, then $F$ has a left $D$-separated subset $L \subseteq F$ with \[|L| \geq |F|/|D|^2.\]
	\end{lemma}
	\begin{scproof}
		First recall some basic facts about finite graphs. Let $G$ be a finite graph with vertex set $V$ and edge set $E$. A subset $I \subseteq V$ is \emphd{independent \ep{in $G$}} if no two vertices in $I$ are adjacent to each other. A straightforward greedy construction shows that if every vertex of $G$ has at most $d$ neighbors, then $V$ can be partitioned into $d+1$ independent sets. In particular, $G$ has an independent set of size at least $|V|/(d+1)$.
		
		Now let $G$ be the graph with vertex set $F$ in which distinct vertices $\sigma$, $\sigma' \in F$ are adjacent if and only if $D\sigma \cap D\sigma' \neq \0$, or, equivalently, $\sigma' \in D^{-1}D\sigma$. Every $\sigma \in F$ has at most $|D^{-1}D| - 1 \leq |D|^2 - 1$ neighbors in $G$ \ep{we are subtracting $1$ to account for the fact that $\sigma$ is not adjacent to itself}. Therefore, $G$ has an independent set $L$ of size at least $|F|/|D|^2$, and that is precisely what we need.
	\end{scproof}

	For a partial map $\phi \colon \G \rightharpoonup k$ with finite domain, let $V_\phi \subseteq k^\G$ denote the clopen set given by
	\[
	V_\phi \,\defeq\, \set{x \in k^\G \,:\, x \supset \phi}.
	\]
	Since the topology on $k^\G$ is generated by the sets of the form $V_\phi$, it is enough to prove Theorem~\ref{theo:synd} with $V = V_\phi$. So, let $D \subset \G$ be a nonempty finite set and let $\phi \colon D \to k$. Let $\ell_0 \in \N$ be such that
	\begin{equation}\label{eq:endgame}
		e \cdot (1 - k^{-|D|})^\ell \cdot |D|^2\ell^2 \,<\, 1 \qquad \text{for all } \ell \geq \ell_0.
	\end{equation}
	\ep{Such $\ell_0$ exists since the left-hand side of \eqref{eq:endgame} approaches $0$ as $\ell$ goes to infinity.} We shall argue that $n \defeq \ell_0|D|^2$ is as desired.
	
	Take any finite set $F \subset \G$ of size at least $n$ and let $L \subseteq F$ be a set of size at least $|F|/|D|^2 \geq \ell_0$ such that $L^{-1}$ is left $D$-separated \ep{such $L$ exists due to Lemma~\ref{lemma:indep} applied with $F^{-1}$ in place of $F$}. Define a $k$-instance $\B$ over $\G$ as follows. For each $\gamma \in \G$, let $B_\gamma$ be the bad $k$-event with domain $D_\gamma \defeq D L^{-1}\gamma$ consisting of all maps $\psi \colon D_\gamma \to k$ such that:
	\[
	\text{for each } \lambda \in L,\ \text{there is } \delta \in D\ \text{with } \psi(\delta\lambda^{-1}\gamma) \,\neq\, \phi(\delta).
	\]
	Notice that since $L^{-1}$ is left $D$-separated, $D_\gamma$ is equal to the \emph{disjoint} union of the sets $D \lambda^{-1}\gamma$, $\lambda \in L$. Let $\B\defeq \set{B_\gamma \,:\, \gamma \in \G}$. We claim that if $x \colon \G\to k$ is a solution to $\B$, then $x$ is trapped in $L \cdot V_\phi$ \ep{and hence also in $F \cdot V_\phi$}. Indeed, consider any $\gamma \in \G$. Since $x$ avoids $B_\gamma$, the following holds:
	\[
		\text{there is } \lambda \in L\ \text{such that for all } \delta \in D, \ x(\delta \lambda^{-1}\gamma) \,=\, \phi(\delta).
	\]
	But this precisely means that for some $\lambda \in L$, $\lambda^{-1}\gamma \cdot x \supset \phi$, i.e., $\gamma \cdot x \in L \cdot V_\phi$. In view of Theorem~\ref{theo:LLL}, it only remains to show that $\B$ is correct for the LLL.
	
	First we bound $d(\B)$. Consider any $\gamma \in \G$. We need to find an upper bound on the size of the set
	\[
		\set{\gamma' \in \G \setminus \set{\gamma} \,:\, D_\gamma \cap D_{\gamma'} \neq \0}.
	\]
	To this end, notice that $D_\gamma \cap D_{\gamma'} \neq \0$ if and only if $\gamma' \in LD^{-1}DL^{-1}\gamma$. Thus,
	\[
		 \left|\set{\gamma' \in \G \setminus \set{\gamma} \,:\, D_\gamma \cap D_{\gamma'} \neq \0}\right| \,=\, |LD^{-1}DL^{-1}| - 1 \,\leq\, |D|^2|L|^2 - 1,
	\]
	and hence we conclude that $d(\B) \leq |D|^2|L|^2 - 1$. Next we bound $p(\B)$. Again, fix $\gamma \in \G$ and draw $f \colon D_\gamma \to k$ uniformly at random. Then for each $\lambda \in L$,
	\[
		\P[f(\delta \lambda^{-1} \gamma) \neq \phi(\delta) \text{ for some } \delta \in D] \,=\, 1 - k^{-|D|},
	\]
	and since the sets $D \lambda^{-1}\gamma$, $\lambda \in L$ are disjoint, we conclude that
	\[
		\P[B_\gamma] \,=\, (1 - k^{-|D|})^{|L|}.
	\]
	Therefore, $p(\B) \leq (1 - k^{-|D|})^{|L|}$. To put everything together, $\B$ is correct as long as
	\[
		e \cdot (1 - k^{-|D|})^{|L|} \cdot |D|^2|L|^2 \,<\, 1.
	\]
	Since $|L| \geq \ell_0$, \eqref{eq:endgame} gives the desired result.
	
	\subsubsection*{Acknowledgments} I am grateful to Todor Tsankov and Andy Zucker for insightful discussions and to Dima Sinapova for providing stimulating and productive environment during the \emph{Logic Fest in the Windy City} conference on May 30--June 2, 2019 at the University of Illinois at Chicago. I am also grateful to the anonymous referee for helpful comments.

	\printbibliography

\end{document}